\documentclass{article}

\usepackage{amsmath,amssymb,graphicx} 
\usepackage[margin=1in]{geometry}
\usepackage{enumerate}
\begin{document}

\nocite{*} 

\title{Directed  and irreversible path in Euclidean spaces}

\author{Khashayar Rahimi }

\date{November 12, 2020} 

\maketitle

\begin{abstract}
The aim of this very short note is to relate the directed paths in ${\stackrel{\rm \longrightarrow}{\rm \mathbb{R}^n}}$ to the irreversible paths in ${\stackrel{\rm ir}{\rm \mathbb{R}^n}}$. We first show that there is a directed path from $x$ to $y$ in ${\stackrel{\rm \longrightarrow}{\rm \mathbb{R}^n}}$ iff there exists an irreversible path with same initial and terminal points in ${\stackrel{\rm ir}{\rm \mathbb{R}^n}}$. Also, we prove that every directed path in ${\stackrel{\rm \longrightarrow}{\rm \mathbb{R}^n}}$ is an irreversible path in ${\stackrel{\rm ir}{\rm \mathbb{R}^n}}$.
\end{abstract}

\section{Introduction}
Now we briefly present some preliminaries including definitions and notations.
\vspace{0.5 cm}

\textbf{Definition 1.} ([1]) A directed topological space, or a d-space $X = (X,PX)$ is a topological space equipped with a set $PX$ of continuous maps $\gamma : I \to X$ (where $I=[0,1]$ equipped with subspace topology of standard topology on $\mathbb{R}$), called directed paths or d-paths, satisfying these axioms:
\vspace{.3 cm}

1. every constant map $I \to X$ is directed\\

2. $PX$ is closed under composition with continuous non-decreasing maps from $I$ to $I$\\

3. $PX$ is closed under concatenation
\vspace{0.5 cm}

\textbf{Notation} ([2]) $\stackrel{\rm ir}{\rm \mathbb{R}}$ is the real line equipped with the left order topology, and $\stackrel{\rm ir}{\rm I}$ is $[0,1]$ equipped with subspace topology of $\stackrel{\rm ir}{\rm \mathbb{R}}$.
\vspace{.5 cm}

\textbf{Definition 2.} Let $X$ be a topological space. A function  $\gamma : \stackrel{\rm ir}{\rm I} \to X$ is called an ir-path in $X$, if it is continuous on $\stackrel{\rm ir}{\rm I}$. Also, $\gamma(0)$ is the initial point and  $\gamma(1)$ is the terminal point of the  ir-path $\gamma$. 
\vspace{.5 cm}

\textbf{Definition 3.} ([3]) ${\stackrel{\rm \longrightarrow}{\rm \Gamma_X}} = \bigg \{ (x,y) \ \bigg | \ \exists \gamma \in PX, \ \gamma(0) = x , \ \gamma(1)=y \bigg \}$
\vspace{.5 cm}

\textbf{Definition 4.} ${\stackrel{\rm ir}{\rm \Gamma_X}} = \bigg \{ (x,y) \ \bigg | \ \exists \gamma \in X^{\stackrel{\rm ir}{\rm I}}, \ \gamma(0) = x , \ \gamma(1)=y \bigg \} $

\section{Relation between d-paths and ir-paths}
\textbf{Proposition 1.} ${\stackrel{\rm \longrightarrow}{\rm \Gamma_{\stackrel{\rm \longrightarrow}{\rm \mathbb{R}^n}}}} = {\stackrel{\rm ir}{\rm \Gamma_{\stackrel{\rm ir}{\rm \mathbb{R}^n}}}} $.\\

\textit{Proof.} From [4] we know that\\
\begin{center}
    ${\stackrel{\rm \longrightarrow}{\rm \Gamma_{\stackrel{\rm \longrightarrow}{\rm \mathbb{R}^n}}}} = \bigg \{(x_1 \ldots x_n, y_1 \ldots y_n) \ \bigg | \ x_i \le y_i \ \forall i \in \{1, \ldots , n \}  \bigg \}$
\end{center}
Also, we know from [2] that there exists an ir-path from $x=(x_1 \ldots x_n)  \in \ \stackrel{\rm ir}{\rm \mathbb{R}^n} $ to $y=(y_1 \ldots y_n) \in \ \stackrel{\rm ir}{\rm \mathbb{R}^n}$ iff $y \in \overline{\{x\}}$. Therefore $(x,y) \in \ {\stackrel{\rm ir}{\rm \Gamma_{\stackrel{\rm ir}{\rm \mathbb{R}^n}}}}$ iff $(y_1 \ldots y_n) \in \ \prod_{i=1}^n [x_i, \infty)$. Clearly, for all $i \in \{1, \ldots , n \}$ we have $x_i \le y_i$, which proves the statement. 
\begin{flushright}
$\square$
\end{flushright}

\textbf{Theorem 1.} Every d-path in $\stackrel{\rm \longrightarrow}{\rm \mathbb{R}^n}$ is an ir-path in ${\stackrel{\rm ir}{\rm \mathbb{R}^n}}$.\\

\textit{Proof.} We know that d-paths in $\stackrel{\rm \longrightarrow}{\rm \mathbb{R}^n}$ are non-decreasing paths in $\mathbb{R}^n$. Thus, it suffices to show that every non-decreasing path in $\mathbb{R}^n$ is an ir-path in ${\stackrel{\rm ir}{\rm \mathbb{R}^n}}$.\\

Suppose that $\gamma: I \to \mathbb{R}^n$ is a non-decreasing path from $x=(x_1 \ldots x_n)$ to $y=(y_1 \ldots y_n)$ and $\prod_{i=1}^n (- \infty, m_i)$ is a base element of ${\stackrel{\rm ir}{\rm \mathbb{R}^n}}$. Now, if ${\gamma}^{-1} \bigg (\prod_{i=1}^n (- \infty, m_i) \bigg )$ is an open subset of $\stackrel{\rm ir}{\rm I }$, then $\gamma$ is an ir-path in ${\stackrel{\rm ir}{\rm \mathbb{R}^n}}$. There are two cases for $m_i$. If for some $1 \le i \le n$, $m_i \le x_i$, then ${\gamma}^{-1} \bigg (\prod_{i=1}^n (- \infty, m_i) \bigg ) = \emptyset$ is open in $\stackrel{\rm ir}{\rm I }$. If for all $1 \le i \le n$, $m_i \ge x_i$, then $0 \in {\gamma}^{-1} \bigg (\prod_{i=1}^n (- \infty, m_i) \bigg ) $. Also we know that ${\gamma}^{-1} \bigg (\prod_{i=1}^n (- \infty, m_i) \bigg ) $ is open in $I$. Hence, ${\gamma}^{-1} \bigg (\prod_{i=1}^n (- \infty, m_i) \bigg ) = [0, n)$ is open in $\stackrel{\rm ir}{\rm I }$.\\
The case ${\gamma}^{-1} \bigg (\prod_{i=1}^n (- \infty, m_i) \bigg ) = [0, n) \cup U \neq [0,p)$ does not happen. Consider an element $u \in U$. Since $\gamma$ is non-decreasing, $\gamma(n) \le \gamma(u) <  (m_1 \ldots m_n) $, then $n \in \ {\gamma}^{-1} \bigg (\prod_{i=1}^n (- \infty, m_i) \bigg ) $, a contradiction.
\begin{flushright}
$\square$
\end{flushright}

\textbf{Example.} This example demonstrates that there exist ir-paths in ${\stackrel{\rm ir}{\rm \mathbb{R}^n}}$ which are not d-paths in $\stackrel{\rm \longrightarrow}{\rm \mathbb{R}^n}$.\\
Let $\gamma : \stackrel{\rm ir}{\rm I} \to {\stackrel{\rm ir}{\rm \mathbb{R}^n}}$ be an ir-path from $x=(x_1 \ldots x_n)$ to $y=(y_1 \ldots y_n)$ defined by\\
\begin{center}
    $ \gamma(t) =
\bigg \{
	\begin{array}{ll}
		x  &  \hspace{1cm}  0 \le t < 1 \\
		y &  \hspace{1cm} t = 1 \\
	\end{array}$
\end{center}

Since ${\gamma}^{-1} \bigg (\prod_{i=1}^n (x_i,y_i + \epsilon) \bigg ) = \{1\}$ is not open in $I$, $\gamma: I \to \mathbb{R}^n$ is not continuous, and hence not a path. Therefore, $\gamma$ is an ir-path that is not a d-path.

\section*{} \label{bibsection}

\bibliographystyle{plain}
\bibliography{template}

\begin{thebibliography}{MMMMM} 
\bibitem{BP} Marco Grandis, \textit{Directed Algebraic Topology, Models of non-reversible worlds}, Cambridge University Press, 2009.

\bibitem{Z} Khashayar Rahimi, \textit{Irreversible homotopy and a notion of irreversible Lusternik–Schnirelmann category}, arXiv preprint arXiv:2010.11217, 2020.

\bibitem{S} E. Goubault, \textit{On directed homotopy equivalences and a notion of directed topological complexity}, arXiv preprint arXiv:1709.05702, 2017.

\bibitem{E} A.Borat and M.Gran, \textit{Directed topological complexity of spheres}, Journal of Applied and Computational Topology, 2020.


\end{thebibliography}

\end{document}